\documentclass[a4paper,12pt]{article}
\begin{document}

\title{Approximation in M\"untz spaces $M_{\Lambda ,p}$ of $L_p$
functions for $1<p<\infty $ and bases.}
\author{Sergey V. Ludkowski}
\date{17 August 2016}
\maketitle

\begin{abstract}
 M\"untz spaces satisfying the M\"untz and gap conditions are considered. A
Fourier approximation of functions in the M\"untz spaces $M_{\Lambda
,p}$ of $L_p$ functions is studied, where $1<p<\infty $. It is
proved that up to an isomorphism and a change of variables these
spaces are contained in Weil-Nagy's class. Moreover, existence of
Schauder bases in the M\"untz spaces $M_{\Lambda ,p}$ is
investigated. \footnote{key words and phrases: Banach space; M\"untz
space;
isomorphism; Schauder basis; Fourier series.  
\par Acknowledgement: the author was partially supported by DFG project
number LU219/10-1
\\ Addresses: Dep. of Mathematics,
Paderborn University,
\par Warburger str. 100, Paderborn D-33095, Germany
\par and Department of Applied Mathematics,
\par Moscow State Technical University MIREA,
\par av. Vernadsky 78, Moscow 119454, Russia
\par sludkowski@mail.ru
}
\end{abstract}

\section{Introduction.} The immense branch of functional analysis
is devoted to topological and geometric properties of topological
vector spaces (see, for example,
\cite{jarchb,ludksmj,ludjms09,naribeckb}). Studies of bases in
Banach spaces compose a large part of it (see, for example,
\cite{jarchb,lindliorb,lusky92}-\cite{lusky04,woytb} and references
therein). It is not surprising that for concrete classes of Banach
spaces many open problems remain, particularly for the M\"untz
spaces $M_{\Lambda ,p}$, where $1<p<\infty $ (see
\cite{ialalam2008}-\cite{clarkerd43}, \cite{gurlusb,schwartzb59} and
references therein). These spaces are defined as completions of the
linear span over $\bf R$ or $\bf C$ of monomials $t^{\lambda }$ with
$\lambda \in \Lambda $ on the segment $[0,1]$ relative to the $L_p$
norm, where $\Lambda \subset [0, \infty )$, $~t\in [0,1]$. In his
classical work K. Weierstrass had proved in 1885 the theorem about
polynomial approximations of continuous functions on the segment.
But the space of continuous functions also forms the algebra.
Generalizations of such spaces were considered by C. M\"untz in 1914
such that his spaces had not the algebra structure. The problem was
whether they have bases. Then a progress was for lacunary M\"untz
spaces satisfying the condition $\underline{\lim}_{n\to \infty }
\lambda _{n+1}/\lambda _n>1$ with a countable set $\Lambda $, but in
its generality this problem was not solved \cite{gurlusb}. It is
worth to mention that the system $ \{ t^{\lambda }: \lambda \in
\Lambda \} $ itself does not contain a Schauder basis for a
nonlacunary set $\Lambda $ satisfying the M\"untz and gap
conditions.
\par In section 2
M\"untz spaces satisfying the M\"untz and gap conditions are
considered. A Fourier approximation of functions in the M\"untz
spaces $M_{\Lambda ,p}$ of $L_p$ functions is studied, where
$1<p<\infty $. Necessary definitions are recalled. It is proved that
up to an isomorphism and a change of variables these spaces are
contained in Weil-Nagy's class. For this purpose in Lemmas 3 and 4,
Theorem 5 and Corollary 9 some isomorphisms of M\"untz spaces are
given. Then in Theorem 13 a relation between M\"untz spaces and
Weil-Nagy's classes is established. Moreover, existence of Schauder
bases in the M\"untz spaces $M_{\Lambda ,p}$ is investigated in
Theorem 16 with the help of Fourier series approximation (see Lemma
14). There is proved that under the M\"untz condition and the gap
condition Schauder bases exist in the M\"untz spaces $M_{\Lambda
,p}$, where $1<p<\infty $.
\par All main results of this paper are obtained for the first time.
They can be used for further investigations of function
approximations and geometry of Banach spaces. It is important not
only for development of mathematical analysis and functional
analysis, but also in their many-sided applications.
\section{Approximation in M\"untz $L_p$ spaces.}
\par To avoid misunderstandings we first remind necessary
definitions and notations.
\par {\bf 1. Notation.} Let $C([\alpha ,\beta ], {\bf F})$ denote the Banach space
of all continuous functions $f: [\alpha ,\beta ]\to \bf F$ supplied
with the absolute maximum norm \par $ \| f \|_C := \max \{ |f(x)|:
x\in [\alpha ,\beta ] \} $, \\ where $- \infty < \alpha < \beta <
\infty $, $~{\bf F}$ is either the real field $\bf R$ or the complex
field $\bf C$.
\par Then $L_p([\alpha ,\beta ],{\bf F})$ denotes the Banach space of
all Lebesgue measurable functions $f: [\alpha , \beta ] \to \bf F$
possessing the finite norm $$ \| f \|_{L_p([\alpha , \beta ], {\bf
F})} := (\int_{\alpha }^{\beta } |f(x)|^pdx)^{1/p} < \infty ,$$
where $1\le p<\infty $ is a marked number, $\alpha <\beta $.
\par Suppose that $Q = (q_{n,k})$ is a lower triangular infinite matrix with real matrix elements
$q_{n,k}$ satisfying the restrictions: $q_{n,k}=0$ for each $k>n$,
where $k, n$ are nonnegative integers. To each $1$-periodic function
$f: {\bf R}\to {\bf R}$ in the space $L_p((\alpha ,\alpha +1),{\bf
F})$ or in $C_0([\alpha ,\alpha +1],{\bf F}) := \{ f: f \in
C([\alpha ,\alpha +1],{\bf F}), f(\alpha )=f(\alpha +1) \} $ is
posed a trigonometric polynomial
$$(1)\quad U_n(f,x,Q) := \frac{a_0}{2} q_{n,0} + \sum_{k=1}^n q_{n,k} (a_k
\cos (2\pi kx) + b_k \sin (2\pi kx)),$$ where $a_k=a_k(f)$ and
$b=b_k(f)$ are the Fourier coefficients of a function $f(x)$.
\par For measurable $1$-periodic functions $h$ and $g$ their convolution is defined whenever it exists by the formula:
$$(2)\quad (h*g)(x) := 2 \int_{\alpha }^{\alpha +1}
h(x-t)g(t)dt.$$
Putting the kernel of the operator $U_n$ to be:
$$(3)\quad U_n(x,Q) := \frac{q_{n,0}}{2}+\sum_{k=1}^n q_{n,k} \cos (2\pi kx)$$
we get
$$(4)\quad U_n(f,x,Q) = (f*U_n(,Q))(x)= (U_n(,Q)*f)(x).$$
The norms of these operators are:
$$(5)\quad L_n(Q,E) := \sup_{f\in E, ~ \| f \| _E=1} \| U_n(f,x,Q) \|_E ,$$
which are constants of a summation method, where $\|*\|_E$ denotes a
norm on a Banach space $E$, where either $E=C_0([\alpha ,\alpha +1],{\bf
F})$ or $E=L_p((\alpha ,\alpha +1),{\bf F})$ with $1\le p<\infty $,
while $\alpha \in {\bf R}$ is a marked real number.
\par As usually $span_{\bf F}(v_k: k)$ will stand for the linear span of
vectors $v_k$ over a field $\bf F$.
\par Henceforward the Fourier summation methods prescribed by sequences of operators
$ \{ U_m : m \} $ which converge on $E$
$$(6)\quad \lim_{m\to \infty } U_m(f,x,Q)=f(x)$$
in the $E$ norm will be considered.

\par {\bf 2. Definition.} Take a countable infinite subset $\Lambda = \{ \lambda_k: k\in{\bf N} \} $
in the set $(0,\infty )$ so that $ \{ \lambda_k: k\in{\bf N} \} $ is
a strictly increasing sequence.
\par Henceforth it is supposed that the set $\Lambda $ satisfies
the gap condition
\par $(1)\quad \inf_k \{ \lambda _{k+1} - \lambda _k \} =: \alpha _0>0$
and the M\"untz condition
\par $$(2)\quad \sum_{k=1}^{\infty } \frac{1}{\lambda _k}=:\alpha _1<\infty .$$
\par The completion of the linear space containing all monomials
$a t^{\lambda }$ with $a\in \bf F$ and $\lambda \in \Lambda $ and
$t\in [\alpha ,\beta ]$ relative to the $L_p$ norm is denoted by
$M_{\Lambda ,p}([\alpha ,\beta ],{\bf F})$, where $0\le \alpha
<\beta <\infty $, $~1\le p$, also by $M_{\Lambda ,C}([\alpha ,\beta
],{\bf F})$ when it is completed relative to the $ \| \| _C$ norm.
Shortly they will also be written as $M_{\Lambda ,p}$ or $M_{\Lambda
,C}$ respectively for $\alpha =0$ and $\beta =1$, when $\bf F$ is
specified.

\par Before subsections about the Fourier approximation in M\"untz
spaces auxiliary Lemmas 3, 4 and Theorem 5 are proved about
isomorphisms of M\"untz spaces $M_{\Lambda ,L_p}$. With the help of
them our consideration reduces to a subclass of M\"untz spaces
$M_{\Lambda ,L_p}$ so that a set $\Lambda $ is contained in the set
of natural numbers $\bf N$.

\par {\bf 3. Lemma.} {\it For each $0<\delta <1$ the M\"untz spaces $M_{\Lambda ,p}([0,1],{\bf F})$ and
$M_{\Lambda ,p}([\delta ,1],{\bf F})$ are linearly topologically
isomorphic, where $1\le p<\infty $.}
\par {\bf Proof.} For every $0<\delta <1$ and $0<\epsilon \le 1$ and $f\in E:=L_p([0,1],{\bf F})$
the norms $\| f \|_{ E[0,1]}$ and $\epsilon \| f|_{[0, \delta ]}
\|_{ E[0,\delta ]} + \| f|_{[\delta , 1]} \|_{ E[\delta , 1]}$ are
equivalent, where $E[\alpha ,\beta ] := E\cap L_p([\alpha ,\beta
],{\bf F})$ for $0\le \alpha <\beta \le 1$. Due to the Remez-type and the
Nikolski-type inequalities (see Theorem 6.2.2 in \cite {borwerdb}
and Theorem 7.4 in \cite{borerdjams1997}) for each $\Lambda $
satisfying the M\"untz condition there is a constant $\eta >0$ so
that $\| h|_{[0, \delta ]} \|_{ E[0,\delta ]} \le \eta \|
h|_{[\delta , 1]} \|_{ E[\delta , 1]}$ for each $h\in M_{\Lambda
,p}$, where $\eta $ is independent of $h$. Therefore the norms $\|
h|_{[\delta , 1]} \|_{ E[\delta , 1]}$ and $ \| h \|_{ E[0, 1]}$ are
equivalent on $M_{\Lambda ,p}[0,1]$. Certainly each polynomial
$a_1t^{\lambda _1}+...+a_nt^{\lambda _n}$ defined on the segment
$[\delta ,1]$ has the natural extension on $[0,1]$, where
$a_1,...,a_n\in {\bf F}$ are constants and $t$ is a variable. Thus
the M\"untz spaces $M_{\Lambda ,p}[0,1]$ and $M_{\Lambda ,p}[\delta ,1]$
are linearly topologically isomorphic as normed spaces for each
$0<\delta <1$.

\par {\bf 4. Lemma.} {\it The M\"untz spaces $M_{\Lambda ,p}$ and
$M_{\Xi \cup (\alpha \Lambda +\beta ),p}$ are linearly topologically
isomorphic for every $ \beta \ge 0$ and $\alpha >0$ and a finite
subset $\Xi $ in $(0,\infty )$, where $1\le p<\infty $.}
\par {\bf Proof.} We have that a sequence $\{ \lambda_k: k\in{\bf N} \} $
is strictly increasing and satisfies the gap condition and hence
$\lim_n \lambda _n = \infty $. We order a set $\Xi \cup (\alpha
\Lambda +\beta )$ into a strictly increasing sequence also. \par  In
virtue of Theorem 9.1.6 \cite{gurlusb} the M\"untz space $M_{\Lambda
, p}$ contains a complemented isomorphic copy of $l_p$,
consequently, $M_{\Lambda ,p}$ and $M_{\Xi \cup \Lambda ,p}$ are
linearly topologically isomorphic as normed spaces.
\par Then from Lemma 3 taking $\alpha
>0$ we deduce that $$(1)\quad \int_{\delta }^1|f(t)|^pdt= \alpha \int_{\delta ^{1/\alpha }}^1|f(x^{\alpha
})|^px^{({\alpha }-1)}dx \le  \alpha \max (1, \delta ^{(1-{\alpha
}^{-1})}) \int_{\delta ^{1/\alpha }}^1|f(x^{\alpha })|^pdx$$ and
$$(2)\quad \int_{\delta }^1|f(x^{\alpha })|^pdx=\alpha ^{-1}\int_{\delta
^{\alpha } }^1|f(t)|^pt^{(\alpha ^{-1}-1)}dt\le \alpha ^{-1}\max
(1,\delta ^{(1- \alpha )})\int_{\delta ^{\alpha } }^1|f(t)|^pdt$$
for each $f\in M_{\Lambda ,p}$ and hence $M_{\alpha \Lambda ,p}$ is
isomorphic with $M_{\Lambda ,p}$. Considering the set $\Lambda _1 =
\Lambda \cup \{ \frac{\beta }{\alpha } \} $ and then the set $\alpha
\Lambda _1$ we get that $M_{\Lambda ,p}$ and $M_{\alpha \Lambda
+\beta ,p}$ are linearly topologically isomorphic as normed spaces
as well.

\par {\bf 5. Theorem.} {\it Let increasing sequences $\Lambda = \{ \lambda
_n : n \} $ and $ \Upsilon = \{ \upsilon _n: n \} $ of positive
numbers satisfy Conditions 2$(1,2)$ and let $\lambda _n\le \upsilon _n$
for each $n$. If $\sup_n (\upsilon _n-\lambda _n)=\delta $, where
$\delta <(8\sum_{n=1}^{\infty }\lambda _n^{-1})^{-1}$, then
$M_{\Lambda ,p}$ and $M_{\Upsilon, p}$ are the isomorphic Banach
spaces, where $1\le p<\infty $.}
\par {\bf Proof.} There exist the natural isometric linear embeddings of
the M\"untz spaces $M_{\Lambda ,p}$ and $M_{\Upsilon, p}$ into $M_{\Lambda \cup
\Upsilon ,p}$. We choose a sequence of sets $\Upsilon _k$ satisfying
the following restrictions \par $(1)$ $\Upsilon _k= \{
\upsilon_{k,n}: n=1,2,... \} \subset \Lambda \cup \Upsilon $ and
$\upsilon _{k,n} \in \{ \lambda _n, \upsilon _n \} $ for each
$k=0,1,2,...$ and $n=1,2,...$, where $\Upsilon _0=\Lambda $;
\par $(2)$ $\upsilon _{k,n}\le \upsilon _{k+1,n}$ for each $k=0,1,2,...$ and
$n=1,2,...$;
\par $(3)$ $\{ \Delta _{k+1,n}: n=1,2,... \} $ is a monotone decreasing
subsequence tending to zero (may be finite or infinite) with
positive terms $\Delta _{k+1,n}$ obtained from the sequence $\delta
_{k+1,j} := \upsilon _{k+1,j}-\upsilon _{k,j}$ by elimination of
zero terms. Denote by $\theta = \theta _{k+1} : \{ j: j\in {\bf N},
\delta _{k+1,j}\ne 0 \} \to {\bf N}$ the corresponding enumeration
mapping such that $\Delta _{k+1,\theta (j)} = \delta _{k+1,j}$ for
each $j\in {\bf N}$ so that $\delta _{k+1,j}\ne 0$ is not zero;
\par $(4)$ $\{ m(k+1): k \} $ is a monotone increasing sequence with
$m(k+1) := \min \{ n: \upsilon _n -\upsilon _{k+1,n}\ne 0; \forall
l<n ~ \upsilon _l =\upsilon _{k+1,l} \} $.
\par Let $f\in M_{\Upsilon _k,p}$. In view of Theorem 6.2.3 and
Corollary 6.2.4 \cite{gurlusb} a function $f$ has a power series
expansion
\par $f(z)=\sum_{n=1}^{\infty } a_nz^{v_{k,n}}$ on $[0,1)$, \\
where $a_n\in {\bf F}$ for each $n\in \bf N$. \par Therefore, for
each $f\in M_{\Upsilon _k, p}$ we consider the power series $f_1(t)
= \sum_{n=1}^{\infty } a_n t^{\upsilon _{k+1,n}}$, where the power
series decomposition $f(t) = \sum_{n=1}^{\infty } a_n t^{\upsilon
_{k,n}}$ converges for each $0\le t <1$, since $f$ is analytic on
$[0,1)$. Then we infer that
$$f(t^2)-f_1(t^2) = \sum_{n=1}^{\infty } a_n t^{\upsilon
_{k,n}}u_{\theta (n)}(t)\mbox{ with }u_{\theta (n)}(t) :=
t^{\upsilon _{k,n}}- t^{\upsilon _{k,n}+2\Delta _{k+1,\theta (n)}}$$
so that $u_l(t)$ is a monotone decreasing sequence by $l$ and hence
\par $|f(t^2)-f_1(t^2)|\le 2 |u_{\theta (m(k+1))}(t)| |f(t)|$ \\
according to Dirichlet's criterium (see, for example, \cite{fichtb})
for each $0\le t<1$, where $\theta =\theta _{k+1}$. Therefore, the
function $f_1(t)$ is analytic on $[0,1)$ and
$$ (5)\quad \|
f-f_1 \| _{L_p([0,1],{\bf F})} \le 2^{2+1/p}  \| f \|
_{L_p([0,1],{\bf F})} \Delta _{k+1,\theta (m(k+1))}/ \lambda
_{m(k+1)},$$ since the mapping $t\mapsto t^2$ is the orientation
preserving diffeomorphism of $[0,1]$ onto itself, also
$|u_{m(k+1)}(t)|\le 2 \Delta _{k+1,\theta (m(k+1))}/ \lambda
_{m(k+1)}$ for each $0\le t\le 1$ by Lemma 7.3.1 \cite{gurlusb} and
$$ \| f-f_1 \| _{L_p([0,1],{\bf F})} = [\int_0^1 |f(\tau ) -
f_1(\tau )|^pd\tau ]^{1/p} $$ $$ = [2\int_0^1 |f(t^2) -
f_1(t^2)|^ptdt]^{1/p} \le [2^{p+1}\int_0^1
|u_{m(k+1)}(t)|^p|f(t)|^ptdt]^{1/p}$$
$$\le 2^{2+1/p}[\int_0^1|f(t)|^pdt]^{1/p} \Delta _{k+1,\theta (m(k+1))}/ \lambda
_{m(k+1)}.$$   Thus the series $\sum_{n=1}^{\infty } a_n t^{\upsilon
_{k+1,n}}$ converges on $[0,1)$. \par Inequality $(5)$ implies that
the linear isomorphism $T_k$ of $M_{\Upsilon _k, p}$ with
$M_{\Upsilon _{k+1}, p}$ exists such that $ \| T_k -I \| \le
2^{2+1/p} \Delta _{k+1,\theta (m(k+1))}/ \lambda _{m(k+1)}$, $ ~
T_k: M_{\Upsilon _k, p}\to M_{\Upsilon _{k+1}, p}$. Then we take the
sequence of operators $S_n := T_nT_{n-1}...T_0: M_{\Lambda , p}\to
M_{\Upsilon _{n+1}, p}\subset M_{\Lambda \cup \Upsilon ,p}$. The
space $M_{\Lambda \cup \Upsilon ,p}$ is complete and the sequence
$\{ S_n: n \} $ operator norm converges to an operator $S:
M_{\Lambda , p}\to M_{\Lambda \cup \Upsilon ,p}$ so that $ \| S-I \|
<1$, since
$$\sum_{k=0}^{\infty }\Delta _{k+1,\theta (m(k+1))}/ \lambda _{m(k+1)}\le
\delta \sum_{n=1}^{\infty }\lambda _n^{-1}<1/8 $$ and $p\ge 1$,
where $I$ denotes the unit operator. Therefore, the operator $S$ is
invertible. On the other hand, from Conditions $(1-4)$ it follows
that $S(M_{\Lambda , p}) =M_{\Upsilon ,p}$.
\par Now we recall necessary definitions and notations of the
Fourier approximation theory and then present useful lemmas.
\par {\bf 6. Notation.} Henceforth $\sf F$ denotes the set of all pairs $(\psi , \beta )$ satisfying
the conditions: $(\psi (k): ~ k\in {\bf N} )$ is a sequence of
non-zero numbers for which $\lim_{k\to \infty }\psi (k)=0$ the limit
is zero, $\beta $ is a real number, also
$$(1)\quad {\cal D}_{\psi , \beta }(x) := \sum_{k=1}^{\infty } \psi (k) \cos (2\pi kx+\beta \pi /2)$$
is the Fourier series of some function from $L_1[0, 1]$. By ${\sf
F}_1$ is denoted the family of all positive sequences $(\psi (k): ~
k\in {\bf N} )$ tending to zero with $\Delta_2\psi (k) := \psi (k-1)
- 2\psi (k) +\psi (k+1)\ge 0$ for each $k$ so that the series
$$(2)\quad \sum_{k=1}^{\infty } \frac{\psi (k)}{k} <\infty $$
converges. The set of all downward convex functions $\psi (v)$ for
each $v\ge 1$ so that $\lim_{v\to \infty } \psi (v)=0$ is denoted by
$\cal M$, while ${\cal M}_1$ is its subset of functions satisfying
Condition $(2)$.
\par Then $$(3)\quad \rho _n (f,x) := f(x) - S_{n-1}(f,x)$$ is the approximation
precision of $f$ by the Fourier series $S(f,x)$, where
$$(4)\quad S_n(f,x) := \frac{a_0}{2} + \sum_{k=1}^n (a_k
\cos (2\pi kx) + b_k \sin (2\pi kx))$$ is the partial Fourier sum
approximating a Lebesgue integrable $1$-periodic function $f$ on
$(0,1)$.

\par {\bf 7. Definition.} Suppose that $f\in L_1(\alpha ,\alpha +1)$ and
$S[f]$ is its Fourier series with coefficients $a_k=a_k(f)$ and
$b_k=b_k(f)$, while $\psi (k)$ is an arbitrary sequence real or
complex. If the function
$$D^{\psi }_{\beta }f := f^{\psi }_{\beta } := \sum_{k=1}^{\infty } [a_k(f)\cos (2\pi
kx+\beta \pi /2) + b_k(f) \sin (2\pi kx + \beta \pi /2)]/\psi (k)$$
belongs to the space $L(\alpha ,\alpha +1)$ of all Lebesgue
integrable (summable) functions on $(\alpha ,\alpha +1)$, then $
f^{\psi }_{\beta }$ is called the Weil $(\psi ,\beta )$ derivative
of $f$. Then $L^{\psi }_{\beta }=L^{\psi }_{\beta }(\alpha ,\alpha
+1)$ stands for the family of all functions $f\in L(\alpha , \alpha
+1)$ with $f^{\psi }_{\beta }\in L(\alpha , \alpha +1)$, we also put
$L^{\psi }_{\beta ,p} := \{ f: f \in L^{\psi }_{\beta }, ~ \|
f^{\psi }_{\beta } \|_{L_p(\alpha ,\alpha +1)}\le 1 \} $.
Particularly, for $\psi (k) =k^{-r}$ this space $L^{\psi }_{\beta }$
is Weil-Nagy's class $W^r_{\beta }=W^r_{\beta }(\alpha , \alpha +1)$
and the notation $W^r_{\beta ,p}$ can be used instead of $L^{\psi
}_{\beta ,p}$ in this case. Put particularly $W^r_{\beta }L_p(\alpha
, \alpha +1) := \{ f: f\in L_p(\alpha , \alpha +1), \exists f^{\psi
}_{\beta }\in L_p(\alpha , \alpha +1) \} $, where $1<p<\infty $.
\par Then let ${\cal E}_n (X) := \sup \{ \| \rho _n(f;x) \|
_{L_p(\alpha ,\alpha +1)}: f \in X \} $,
\par $E_n(f)_p := \inf \{ \| f- T_{n-1} \|_{L_p(\alpha ,\alpha +1)}: T_{n-1}\in
{\cal T}_{2n-1} \} $,
\par $E_n(X) := \sup \{ E_n(f)_p: f \in X
\} $, \\ where $X$ is a subset in $L_p(\alpha ,\alpha +1)= L_p((\alpha
,\alpha +1),{\bf R})$,
$${\cal T}_{2n-1} := \{ T_{n-1}(x) = \frac{c_0}{2} +
\sum_{k=1}^{n-1} (c_k \cos (2\pi kx) + d_k \sin (2\pi kx)); ~ c_k,
d_k \in {\bf R} \} $$ denotes the family of all trigonometric
polynomials $T_{n-1}$ of degree not greater than $n-1$.

\par {\bf 8. Lemma.} {\it Suppose that
$Q_{\alpha }f(t) := f(t^{\alpha })$ for each $f: [0,1]\to \bf F$,
where $0<\alpha $, $t\in [0,1]$, $1<p<\infty $. Then for each
$1<\alpha <\infty $ there exists $0<\delta < 1$ such that the
operator $Q_{\alpha }$ from $L_ p(\delta ^{\alpha },1)$ into
$L_p(\delta ,1)$ has the norm $ \| Q_{\alpha } \| <1$.}
\par {\bf Proof.} The Banach spaces $L_p(\delta ,1)$ and
$L_p(\delta ^{\alpha },1)$ are defined with the help of the Lebesgue
measure on $\bf R$. Then Formula 4$(2)$ implies that $ \| Q_{\alpha
} \| <1$ as soon as $\alpha ^{-1} \max (1,\delta ^{(1-\alpha )})<1$.
That is when $\{ \delta ^{(1-\alpha )}< \alpha \}
\Longleftrightarrow \{ \ln \delta > (1- \alpha )^{-1} \ln \alpha \}
$ , since $\alpha >1$ and $0<\delta <1$.

\par {\bf 9. Corollary.} {\it Let $1<\alpha <\infty $ and $0<\delta <1$
so that $\delta > \alpha ^{1/(1- \alpha )}$, let also $Z_{\Lambda ,p
,\alpha , \delta } := (I-Q_{\alpha }) [M_{\Lambda ,p}(\delta
^{\alpha },1)]$, where $1<p<\infty $, while $I$ is the unit
operator. Then $Z_{\Lambda , p, \alpha ,\delta }$ is isomorphic with
$M_{\Lambda ,p}(\delta ^{\alpha },1)$.}
\par {\bf Proof.} There is the natural embedding of $L_p(a,b)$ into
$L_p(c,d)$ when $c\le a$ and $b\le  d$ such that $f\mapsto f\chi
_{(a,b)}$ for each $f \in L_p(a,b)$, where $\chi _A$ notates the
characteristic function of a set $A$. Since $ \| Q_{\alpha } \| <1$,
then the operator $I-Q_{\alpha }$ is invertible (see
\cite{kolmfomb}).

\par {\bf 10. Lemma.} {\it Let $f\in L_p(0,1)$, where $1<p<\infty $. Then $$\lim_{\eta
\downarrow 0} \eta ^{-1/q} \int_{1-\eta }^1 f(t)dt=0,$$ where
$1/q+1/p=1$.}
\par {\bf Proof.} Since $f\in L_p(0,1)$, then $|f(t)|^p\mu (dt)$ is a $\sigma $-additive
and finite measure on $(0,1)$, where $\mu $ is the Lebesgue measure
on $\bf R$ (see, for example, \cite{edwardsb}, Theorems V.5.4.3 and
V.5.4.5 \cite{kolmfomb}). Therefore, the limit exists
$$(1)\quad \lim _{\eta
\downarrow 0} \int_{1-\eta }^1 |f(t)|^pdt=0.$$ From Holder's
inequality it follows that $$|\int_{1-\eta }^1 f(t)dt|\le
(\int_{1-\eta }^1 |f(t)|^pdt)^{1/p} (\int_{1-\eta }^11dt)^{1/q}$$ $$
= \eta ^{1/q}(\int_{1-\eta }^1 |f(t)|^pdt)^{1/p}\mbox{ hence }$$
$$(2)\quad |\eta ^{-1/q}\int_{1-\eta }^1 f(t)dt|\le (\int_{1-\eta }^1
|f(t)|^pdt)^{1/p} .$$ Thus from Formulas $(1)$ and $(2)$ the
statement of this lemma follows.

\par {\bf 11. Note.}  We remind the following
definition: the family of all Lebesgue measurable functions $f:
(a,b)\to \bf R$ satisfying the condition
$$\| f \|_{L_{s,w}(a,b)} :=  \sup_{y>0}(y^s\mu \{ t: ~ t\in (a,b), ~ |f(t)|\ge y \} )^{1/s} <\infty
$$ is called the weak $L_s$ space and denoted by $L_{s,w}(a,b)$, where $\mu $ notates the Lebesgue
measure on the real field $\bf R$, $~0<s<\infty $, $~(a,b)\subset
{\bf R}$ (see, for example, \S 9.5 \cite{edwardsb}, \S IX.4
\cite{reedsim}, \cite{steinb}).
\par The following proposition 12 is used below in theorem 13 to prove that
functions of M\"untz spaces $M_{\Lambda , p}$ for $\Lambda $
satisfying the M\"untz condition and the gap condition belong to
Weil-Nagy's class, where $1<p<\infty $.

\par {\bf 12. Proposition.} {\it Suppose that an increasing sequence
$\Lambda = \{ \lambda _n : n \} $ of natural numbers satisfies the
M\"untz condition, $1<p<\infty $ and $f\in M_{\Lambda ,p}$. Then
$dh(x)/dx\in L_{s,w}(0,1)$ for a function $h(x)=f(x)-f(x^2)$, where
$s= p/(p+1)$.}

\par {\bf Proof.} By the conditions of this proposition $f\in M_{\Lambda
,p}$. By virtue of \cite{clarkerd43} a function $f$ has the analytic
extension on $\dot{B}_1(0)$ and the series
\par $(1)$ $f(z)=\sum_{n=1}^{\infty } a_nz^{\lambda _n}$\\
converges on $\dot{B}_1(0)$, where $\dot{B}_r(x) := \{ y: y\in {\bf
C}, |y-x|<r \} $ denotes the open disk in $\bf C$ of radius $r>0$
with center at $x\in {\bf C}$, where $a_n\in {\bf F}$ is an
expansion coefficient for each $n\in \bf N$. This also evidently
follows from Theorems 6.2.2, 6.2.3 and Corollary 6.2.4
\cite{gurlusb}, the Abel theorem 3 about a power series on a disk
and the Cauchy-Hadamard formula of a power series convergence radius
and Theorem 5 in Subsection 20, Section 6, Chapter II in
\cite{shabatb}, since $\Lambda \subset {\bf N}$. That is, the
function $f$ has a holomorphic univalent extension from $[0,1)$ on
$\dot{B}_1(0)$. \par Take the function $H(x) =\int_x^1h(t)dt$, where
$x\in [0,1]$. In view of Theorem VI.4.2 \cite{kolmfomb} and
Lyapunov's inequality (formula $(27)$ in \S II.6 \cite{shirb}) this
function is continuous so that $H(1)=0$. Together with formula $(1)$
this implies that the function $H(x)$ belongs to $M_{\{ 0 \} \cup
(\Lambda +1),C}$ and has a holomorphic univalent extension on
$\dot{B}_1(0)$.
\par Then we put $g(z) = (1-z)^{-1/q} H(z)$ for each $z\in \dot{B}_1(0)$, where $1/q+1/p=1$.
From Lemma 10 it follows that \par $(2)$ $\lim_{z\to 1} g(z)=0$. \\
Thus the function $g(z)$ is holomorphic (may be multivalent because
of the multiplier $(1-z)^{-1/q}$) on $\dot{B}_1(0)$ and continuous
on $\dot{B}_1(0)\cup \{ 1 \} $.
\par According to Cauchy's formula 21$(5)$ in \cite{shabatb}
$$(3)\quad h'(z) = - \frac{1}{\pi i} \int_{\gamma }
\frac{H(y)}{(y-z)^3}dy$$ for each $z\in \dot{B}_{1/2}(1/2)$, where
$\gamma $ is an oriented rectifiable boundary $\gamma =
\partial G$ of a simply connected open domain
$G$ contained in $\dot{B}_{1/2}(1/2)$ such that $z\in G$.
Particularly, this is valid for each $z$ in $(1/2,1)$ and $G=
\dot{B}_{1/2}(1/2)$. \par On the other hand, the function $g(z)$ is
bounded on $B_{1/2}(1/2)$, where $B_r(x) := \{ y: y\in {\bf C},
|y-x|\le r \} $ notates the closed disk of radius $r>0$ with center
at $x\in {\bf C}$. Thus $K=\sup_{z\in B_{1/2}(1/2)} |g(z)|<\infty $.
Estimating the integral $(3)$ and taking into account formula $(2)$
we infer that $|h'(t)|\le 2K/(1-t)^{1+1/p}$ for each $t\in (3/4,1)$,
since $1/q+1/p=1$. Together with the analyticity of $h'$ on $[0,1)$
this implies that
\par $\sup_{y>0}(y^s\mu \{ t: ~ t\in (a,b), ~ |h'(t)|\ge y \} )^{1/s}
<\infty $, \\ where $s= p/(p+1)$. Thus $h'\in L_{s,w}(0,1)$.

\par {\bf 13. Theorem.} {\it Let an increasing sequence $\Lambda = \{
\lambda _n : n \} $ of natural numbers satisfy the M\"untz
condition, also $1>\delta > 1/2$ and $1<p<\infty $ and let $\sigma
(x) = \delta ^2+ x (1-\delta ^2)$, where $0\le x\le 1$. Then for
each $0<\gamma <1$ there exists $\beta = \beta (\gamma )\in {\bf R}$
so that $Z_{\Lambda ,p ,2, \delta }\circ \sigma \subset W^{\gamma
}_{\beta }L_p(0,1)$.}
\par {\bf Proof.} Let $f\in M_{\Lambda ,p}(0,1)$ and
$v(x)= (I-Q_2)f(\sigma (x))$, then $v(x)$ is analytic on $(0,1)$,
since $f$ is analytic on $(0,1)$ and $\sigma [0,1]=[\delta ^2, 1]$.
We take its $1$-periodic extension $v_0$ on $\bf R$.
\par According to Proposition 1.7.2 \cite{stepanetsb} (or see
\cite{zygmb}) $h\in W^{\gamma }_{\beta }L_p(0,1)$ if and only if
there exists a function $\phi = \phi_{h, \gamma ,\beta }$ which is
$1$-periodic on $\bf R$ and Lebesgue integrable on $[0,1]$ such that
$$(1)\quad h(x) = \frac{a_0(h)}{2} + (\phi * {\cal D}_{\psi ,\beta
})(x),$$ where $a_0(h)= 2 \int_0^1 h(t)dt$ (see \S \S 6 and 7).
\par We take a sequence $U_n(t,Q)$ given by Formula 1$(3)$ so that
$$\lim_m q_{m,k} =1 \mbox { for each }k \mbox{ and } \sup_m
L_m(Q,L_p) < \infty \mbox{ and }  \sup_{m,k} |q_{m,k}| <\infty $$
and write for short $U_n(t)$ instead of $U_n(t,Q)$. Under these
conditions the limit exists
\par $(2)$ $\lim_n (v*U_n)(x) =v(x)$ \\ in $L_p(0,1)$ norm for each
$v\in L_p((0,1),{\bf F})$ according to Chapters 2 and 3 in
\cite{stepanetsb} (see also \cite{barib,zygmb}).
\par On the other hand, Formula I$(10.1)$ \cite{stepanetsb} provides
\par $(3)$ $S[(y^{\psi _1}_{\bar{\beta }_1})^{\psi _2/\psi _1}_{\bar{\beta
}_2 - \bar{\beta }_1}]=S[y^{\psi _2}_{\bar{\beta}_2}]$, \\ where
$S[y]$ is the Fourier series corresponding to a function $y\in
L^{\psi _2}_{\bar{\beta}_2}$, when $(\psi _1,\bar{\beta }_1)\le
(\psi _2,\bar{\beta }_2)$.
\par Put $\theta (k)=k^{\gamma -1}$ for all $k\in \bf N$. Then
${\cal D}_{\theta , - \beta }\in L_1(0,1)$ for each $\beta \in \bf
R$ due to Theorems II.13.7, V.1.5 and V.2.24 \cite{zygmb} (or see
\cite{barib}). This is also seen from chapters I and V in
\cite{stepanetsb} and Formulas $(1)$ and $(3)$ above. In view of
Dirichlet's theorem (see \S 430 in \cite{fichtb}) the function
${\cal D}_{\theta , -\beta }(x) $ is continuous on the segment
$[\delta ,1-\delta ]$ for each $0<\delta <1/4$.
\par According to formula 2.5.3.$(10)$
in \cite{prbrmarb1}
$$\int_0^{\infty } x^{\alpha -1} {\sin (bx)\choose{\cos (bx)}} dx =
b^{-\alpha }\Gamma (\alpha )
 {\sin (\pi \alpha /2) \choose{\cos (\pi \alpha /2)}}$$
 for each $b>0$ and $0<Re (\alpha )<1$. On the other hand,
the integration by parts gives: $$\int_a^{\infty } x^{\alpha -1}
{\sin (bx)\choose{\cos (bx)}} dx = b^{-1}a^{\alpha -1} {\cos
(ab)\choose{- \sin (ab)}} - b^{-1}(\alpha -1)\int_a^{\infty }
x^{\alpha -2} {- \cos (bx)\choose{\sin (bx)}} dx$$ for every $a>0$,
$b>0$ and $0<Re (\alpha )<1$. From formulas V$(2.1)$, theorems
V.2.22 and V.2.24 in \cite{zygmb} (see also \cite{bruijnb,olverb})
we infer the asymptotic expansions
$$\sum_{n=1}^{\infty } n^{-\alpha }\sin (2\pi nx) \approx  (2\pi x)^{\alpha -1} \Gamma (1-\alpha )
\cos (\pi \alpha /2) + \mu x^{\alpha },$$
$$\sum_{n=1}^{\infty } n^{-\alpha }\cos (2\pi nx) \approx  (2\pi x)^{\alpha -1} \Gamma (1-\alpha )
\sin (\pi \alpha /2)+ \nu x^{\alpha }$$ in a small neighborhood
$0<x<\delta $ of zero, where $0<\delta <1/4$, $0<\alpha <1$, $\mu $
and $\nu $ are real constants. Taking $\beta =\alpha = 1-\gamma $ we
get that ${\cal D}_{\theta , -\beta }(x) \in L_{\infty }(0,1).$
\par Evidently, for Lebesgue measurable
functions $f: {\bf R}\to {\bf R}$ and $g: {\bf R}\to {\bf R}$ there
is the equality $\int_{-\infty }^{\infty } f(x-t)\chi _{[0,\infty
)}(x-t)g(t)\chi _{[0,\infty )}(t) dt = \int_0^x f(x-t)g(t)dt$ for
each $x>0$ whenever this integral exists, where $\chi _A$ denotes
the characteristic function of a subset $A$ in $\bf R$ such that
$\chi _A(y)=1$ for each $y\in A$, also $\chi _A(y)=0$ for each $y$
outside $A$, $y\in {\bf R}\setminus A$. Particularly, if $0<x\le T$,
where $0<T<\infty $ is a constant, then $\int_0^x f(x-t)g(t)dt=
\int_0^{\infty } f(x-t)\chi _{[0,T]}(x-t)g(t)\chi _{[0,T]}(t)dt$
(see also \cite{fichtb,kolmfomb}). This is applicable to formula
1$(2)$ putting $\alpha =0$ there and with the help of the equality
\par $\int_0^1 f(x-t)g(t)dt = \int_0^xf(x-t)g(t)dt + \int_0^{1-x}
f_1((1-x)-v)g_1(v)dv$ \\ for each $0\le x\le 1$ and $1$-periodic
functions $f$ and $g$ and using also that $ \| f|_{[a,b]} \| \le \|
f|_{[0,1]} \| = \| f_1|_{[0,1]} \|$ for the considered here types of
norms for each $[a,b]\subset [0,1]$, where $f_1(t)=f(-t)$ and
$g_1(t)=g(-t)$ for each $t\in \bf R$, since \par $\int_x^1
f(x-t)g(t)dt = \int_0^{1-x} f(v-1+x) g(1-v)dv$.
\par Mention that according to the weak Young inequality
\par $(4)$ $ \| \xi *\eta  \|_p \le K_{r,s}
\| \xi \|_r \| \eta  \|_{s,w}$ \\
for each $\xi \in L_r$ and $\eta \in L_{s,w}$, where $1\le p,  r \le
\infty $, $0<s<\infty $ and $r^{-1} + s^{-1} =1 + p^{-1}$,
$~K_{r,s}>0$ is a constant independent of $\xi $ and $\eta $ (see
theorem 9.5.1 in \cite{edwardsb}, \S IX.4 in \cite{reedsim}).
\par In virtue of formula $(3)$, the weak Young
inequality $(4)$ and Proposition 12 there exists a function $s$ in
$L_p(0,1)$ so that
$$s(x)= \lim_n (({\cal D}_{\theta , - \beta
}*U_n)*v_0')(x),$$ where $\beta = 1-\gamma .$ Therefore $\phi_{v_0,
\gamma ,\beta }=s$ and $D^{\psi }_{\beta }v_0=s$ according to $(1)$
and $(3)$. Thus $v_0\in W^{\gamma }_{\beta }L_p(0,1)$.
\par Below Lemma 14 and Proposition 15 are given. They are used in subsection 16 for
proving existence of a Schauder basis. On the other hand, Theorem 13
is utilized that to prove Lemma 14.
\par {\bf 14. Lemma.} {\it  If an increasing sequence $\Lambda $
of natural numbers satisfies the M\"untz condition, also $0<\gamma
<1$ and $1< p<\infty $, $1>\delta > 1/2$,
\par $X= \{ h: ~ h=f\circ \sigma , f\in Z_{\Lambda ,p ,2, \delta };
\| f \|_{L_p((\delta ^2,1),{\bf R})}\le 1 \} $,
\par then a positive constant $\omega  = \omega (p,\gamma )$ exists so that
$$(1)\quad E_n(X)\le {\cal E}_n(X) \le \omega n^{-\gamma }
$$ for each natural number $n\in \bf N$.}
\par {\bf Proof.} Due to Theorem 13 the inclusion is valid $h(x)\in
W^{\gamma }_{\beta }L_p(0,1)$ for each $h\in Z_{\Lambda ,p ,2,
\delta }\circ \sigma $, where $\psi $ is in ${\sf F}_1$ so that
$\psi (k)= k^{-\gamma }$ for each $k\in \bf N$, $\beta =1-\gamma $.
Then $\| h \|_{L_p((0,1),{\bf R})}=(1-\delta ^2)^{-1/p} \| f
\|_{L_p((\delta ^2,1),{\bf R})} \le (1-\delta ^2)^{-1/p} $ for each
$h\in X$, since
$$(2)\quad \int_0^1|h(x)|^pdx = (1-\delta ^2)^{-1}
\int_{\delta ^2}^1|f(t)|^pdt .$$ Therefore, $X\subset (1-\delta
^2)^{-1/p}W^{\gamma }_{\beta }L_p(0,1)$ (see also \S 7), where $bY
:= \{ f: f=bg, g\in Y \} $ for a linear space $Y$ over $\bf R$ and a
marked real number $b$.
\par Then estimate $(1)$ follows from Theorem V.5.3 in \cite{stepanetsb}.

\par {\bf 15. Proposition.} {\it Let $X$ be a Banach space over $\bf R$
and let $Y$ be its Banach subspace so that they fulfill conditions
$(1-4)$ below:
\par $(1)$ there is a sequence $(e_i: i \in {\bf N})$ in $X$
such that $e_1,...,e_n$ are linearly independent vectors and $ \|
e_n \|_X =1$ for each $n$ and
\par $(2)$ there exists a Schauder basis $(z_n: n\in {\bf N})$
in $X$ such that \par $z_n = \sum_{k=1}^n b_{k,n}e_k$ for each $n\in
\bf N$, where $b_{k,n}$ are real coefficients;
\par $(3)$ for every $x\in Y$ and $n\in \bf N$ there exist
$x_1,...,x_n\in {\bf R}$ so that
\par $ \| x - \sum_{i=1}^n x_ie_i \|_X \le s(n) \| x \| $, \\
where $s(n)$ is a strictly monotone decreasing positive function
with \par $\lim_{n\to \infty } s(n)=0$ and
\par $(4)$ $u_n =\sum_{l=m(n)}^{k(n)} u_{n,l} e_l$, \\
where $u_{n,l}\in \bf R$ for each natural numbers $k$ and $l$, where
a sequence $(u_n: n \in {\bf N})$ of normalized vectors in $Y$ is
such that its real linear span is everywhere dense in $Y$ and $1\le
m(n)\le k(n)<\infty $ and $m(n)<m(n+1)$ for each $n\in \bf N$.
\par Then $Y$ has a Schauder basis.}
\par {\bf Proof.} Without loss of generality one can select and enumerate
\par $(5)$ vectors $u_1$,...,$u_n$ so that they are linearly independent in $Y$
for each natural number $n$. By virtue of Theorem $(8.4.8)$ in
\cite{naribeckb} their real linear span $span_{\bf R}
(u_1,....,u_n)$ is complemented in $Y$ for each $n\in \bf N$. Put
$L_{n,\infty } := cl_X span_{\bf R} (u_k: k\ge n)$ and $L_{n,m}:=
cl_X span_{\bf R} (u_k: n\le k\le m)$, where $cl_XA$ denotes the
closure of a subset $A$ in $X$, where $span_{\bf R}A$ denotes the
real linear span of $A$. Since $Y$ is a Banach space and $u_k\in Y$
for each $k$, then $L_{n,\infty } \subset Y$ and $L_{n,m}\subset Y$
for each natural numbers $n$ and $m$. Then we infer that \par
$L_{n,j} \subset span_{\bf R} (e_l: m(n)\le l \le k_{n,j} )$, where
$k_{n,j} := \max (k(l): n\le l \le j)$. \par Take arbitrary vectors
$f\in L_{1,j}$ and $g\in L_{j+1,q}$, where $1\le j<q$. Therefore,
there are real coefficients $f_i$ and $g_i$ such that
\par $ f=\sum_{i=1}^{k_{1,j}} f_ie_i $ and
 \par $g=\sum_{i=m(j+1)}^{k_{j+1,q}} g_ie_i$. Hence due to condition
 $(2)$:
\par $ \| f - \sum_{i=1}^{m(j)} f_ie_i \| _X \le s(m(j)) \| f \| $
and \par $ \| g - \sum_{i=k_{1,j}+1}^{k_{j+1,q}} g_ie_i \|_X \le
s(k_{1,j}+1) \| g \|_X$. \par On the other hand,
\par $ f=  \sum_{i=1}^{m(j)} f_ie_i + \sum_{i=m(j)+1}^{k_{1,j}} f_ie_i $, consequently,
\par $ \| f^{[j+1]} \| \le s(m(j+1))) \| f \| $, where
\par $f^{[j+1]} := \sum_{i=m(j+1)}^{k_{1,j}} f_ie_i$ and  $\sum_{i=a}^bf_ie_i := 0$, when
$a>b$.  \par When $0<\delta <1/4$ and $s(m(j)+1)< \delta $ we infer
using the triangle inequality that $ \| f^{[j+1]} - h \|_X \le
\delta \| f^{[j+1]} \|_X /(1-\delta ) \le \delta s(m(j+1)-1) \| f
\|_X /(1-\delta )$ for the best approximation $h$ of $f^{[j+1]}$ in
$L_{j+1,\infty }$, since $m(j)<m(j+1)$ for each $j$. Therefore, the
inequality $ \| f - g \|_X \ge \| f-f^{[j+1} \|_X -  \| f^{[j+1]} -
g \|_X$ and $s(n)\downarrow 0$ imply that there exists $n_0$ such
that the inclination of $L_{1,j}$ to $L_{j+1,\infty }$ is not less
than $1/2$ for each $j\ge n_0$. Condition $(4)$ implies that
$L_{1,n_0}$ is complemented in $Y$. In virtue of Theorem 1.2.3
\cite{gurlusb} a Schauder basis exists in $Y$.

\par {\bf 16. Theorem.} {\it If a set $\Lambda $ satisfies the M\"untz and gap conditions
and $1<p<\infty $, then the M\"untz space $M_{\Lambda ,p}([0,1],{\bf
F})$ has a Schauder basis.}
\par {\bf Proof.} In view of Lemma 4 and Theorem 5 there is sufficient to prove an existence
of a Schauder basis in the M\"untz space $M_{\Lambda ,p}$ for $\Lambda
\subset {\bf N}$. Mention that if the M\"untz space $M_{\Lambda
,p}([0,1],{\bf R})$ over the real field has the Schauder basis then
$M_{\Lambda ,p}([0,1],{\bf C})$ over the complex field has it as
well. Thus it is sufficient to consider the real field ${\bf F}={\bf R}$.
\par  Let $U_m(x,Q)$ be kernels of the Fourier summation method
in $L_p(0,1)$ as in \S 2.1 such that $$(1)\quad \lim_m q_{m,k} =1
\mbox { for each }k \mbox{ and } \sup_m L_m(Q,L_p) < \infty \mbox{
and } \sup_{m,k} |q_{m,k}| <\infty .$$ For example, Cesaro's
summation method of order 1 can be taken to which Fej\'er kernels
$F_n$ correspond so that the limit \par $\lim_{n\to \infty }F_n * f=f$ \\
converges in $L_p(0,1)$ (see Theorem 19.1 and Corollary 19.2 in
\cite{zaanb}). That is, there exists a Schauder basis $z_n$ in
$L_p(0,1)$ such that
\par $z_{2n}(t) = a_{0,2n} + [\sum_{k=1}^{n-1} (a_{k,2n} \cos (2\pi kt) +
b_{k,2n}\sin (2\pi kt)] + a_{n,2n} \cos (2\pi nt)$ and \par
$z_{2n+1}(t) = a_{0,2n+1} + \sum_{k=1}^n (a_{k,2n+1} \cos (2\pi kt)
+ b_{k,2n+1}\sin (2\pi kt))$ \\ for every $t\in (0,1)$ and $n\in
{\bf N}$, where $a_{k,j}$ and $b_{k,j}$ are real expansion
coefficients. \par In virtue of Theorem 6.2.3 and Corollary 6.2.4
\cite{gurlusb} each function $g\in M_{\Lambda ,p}[0,1]$ has an
analytic extension on $\dot{B}_1(0)$ and hence
$$(2)\quad g(z) = \sum_{n=1}^{\infty } c_n z^{\lambda _n} =
\sum_{k=1}^{\infty } p_k u_k(z)$$ are the convergent series on the
unit open disk $\dot{B}_1(0)$ in $\bf C$ with center at zero (see \S
12), where $\Lambda \subset \bf N$ and $c_n=c_n(g)\in \bf N$,
$~p_n=p_n(g)=c_1+...+c_n$, $~ u_1(z) := z^{\lambda _1}$, $~
u_{n+1}(z) := z^{\lambda _{n+1}} - z^{\lambda _n}$ for each
$n=1,2,...$. On the other hand, the M\"untz spaces $M_{\Lambda
,p}[0,1]$ and $M_{\Lambda ,p}[\delta ^2, 1]$ are isomorphic for each
$0<\delta <1$ (see Lemma 3 above). Therefore we consider
henceforward the M\"untz space $M_{\Lambda ,p}$ on the segment
$[\delta ^2,1]$, where $1>\delta
> 1/2$. Mention that $M_{\Lambda ,p}[\delta
^2, 1]$ and $M_{\Lambda ,p}\circ \sigma [0,1]$ are isomorphic (see
\S 13). Then $Z_{\Lambda ,p ,2, \delta }$ and $Z_{\Lambda ,p ,2,
\delta }\circ \sigma |_{[0,1]}$ are isomorphic as well. In view of
Corollary 9 it is sufficient to prove the existence of a Schauder
basis in $Z_{\Lambda ,p ,2, \delta }\circ \sigma |_{[0,1]}$.
\par Take the finite dimensional subspace $X_n := span_{\bf
R}(u_1,...,u_n)$ in $M_{\Lambda ,p}$, where $n\in \bf N$. Due to
Lemma 4 the Banach space $M_{\lambda ,p}\ominus X_n$ exists and is
isomorphic with $M_{\lambda ,p}$. In virtue of Formula I$(10.1)$
\cite{stepanetsb} $S[(y^{\psi _1}_{\bar{\beta }_1})^{\psi _2/\psi
_1}_{\bar{\beta }_2 - \bar{\beta }_1}]=S[y^{\psi
_2}_{\bar{\beta}_2}]$, where $y\in L^{\psi _2}_{\bar{\beta}_2}$,
when $(\psi _1,\bar{\beta }_1)\le (\psi _2,\bar{\beta }_2)$. \par
Consider the trigonometric polynomials $U_m(f,x,Q)$ for $f\in
(Z_{\Lambda ,p ,2, \delta }\ominus (I-Q_2)X_n)\circ \sigma $, where
$m=1, 2,...$. Put $Y_{K,n}$ to be the $L_p$ completion of the linear
span $span_{\bf R} (U_m(f,x,Q): (m,f)\in K)$, where $K\subset {\bf
N}\times (Z_{\Lambda ,p ,2, \delta }\ominus (I-Q_2)X_n)\circ \sigma
$, $m\in \bf N$, $f\in (Z_{\Lambda ,p ,2, \delta }\ominus
(I-Q_2)X_n)\circ \sigma $. \par It is known (see Proposition 1.7.1
\cite{stepanetsb}) that $f\in L^{\psi }_{\beta }(\alpha ,\alpha +1)$
if and only if there exists $g \in L(\alpha ,\alpha +1)$ so that
$f=\frac{a_0(f)}{2}+{\cal D}_{\psi ,\beta }*g$, where the function
${\cal D}_{\psi ,\beta }$ is prescribed by Formula 6$(1)$, the
constant $a_0(f)$ is as above. In view of Lemma 4 it is sufficient
to consider the case $a_0(f)=0$.
\par There exists a countable subset $\{ f_n: n\in {\bf N} \} $ in
$Z_{\Lambda ,p ,2, \delta }$ such that $f_n\circ \sigma ={\cal
D}_{\psi ,\beta }*g_n$ with $g_n\in L(0,1)$ for each $n\in \bf N$
and so that $span_{\bf R}\{ f_n: n\in {\bf N} \}$ is dense in
$Z_{\Lambda ,p ,\alpha , \delta }$, since $Z_{\Lambda ,p ,2, \delta
}$ is separable. Using Formulas $(1,2)$, Proposition 12 and Lemma 14
we deduce that a countable set $K$ and a sufficiently large natural
number $n_0$ exist so that the Banach space $Y_{K,n_0}$ is
isomorphic with $(Z_{\Lambda ,p ,2, \delta }\ominus (I-Q_2)X_{n_0})$
and $Y_{K,n_0}|_{(0,1)}\subset W^{\gamma }_{\beta }L_p(0,1)$, where
$0<\gamma <1$ and $\beta =1-\gamma $. Therefore, by the construction
above the Banach space $Y_{K,n_0}$ is the $L_p$ completion of the
real linear span of a countable family $(s_l: l\in {\bf N})$ of
trigonometric polynomials $s_l$.
\par Without loss of generality this family can be refined by
induction such that $s_l$ is linearly independent of
$s_1,...,s_{l-1}$ over $\bf F$ for each $l\in \bf N$. With the help
of transpositions in the sequence $\{ s_l: l \in {\bf N} \} $, the
normalization and the Gaussian exclusion algorithm we construct a
sequence $ \{ r_l: l\in {\bf N} \} $ of trigonometric polynomials
which are finite real linear combinations of the initial
trigonometric polynomials $\{ s_l : l\in {\bf N} \} $ and satisfying
the conditions
\par $(3)$ $ \| r_l \|_{L_p(0,1)} =1$ for each $l$;
\par $(4)$ the infinite matrix having $l$-th row of the form
$...,a_{l,k}, b_{l,k}, a_{l,k+1}, b_{l,k+1},...$ for each $l\in \bf
N$ is upper trapezoidal (step), where $$r_l(x) = \frac{a_{l,0}}{2} +
\sum_{k=m(l)}^{n(l)} [a_{l,k} \cos (2\pi kx) + b_{l,k} \sin (2\pi
kx) ]$$ with $a_{l,m(l)}^2+b_{l,m(l)}^2>0$ and
$a_{l,n(l)}^2+b_{l,n(l)}^2>0$, where $1\le m(l)\le n(l)$, $deg
(r_l)=n(l)$, or $r_1(x)= \frac{a_{1,0}}{2}$ when $deg (r_1)=0$;
$a_{l,k}, b_{l,k}\in {\bf R}$ for each $l\in {\bf N}$ and $0\le k\in
{\bf Z}$.
\par Then as $X$ and $Y$ in Proposition 15 we take $X=L_p[0,1]$ and
$Y=Y_{K,n_0}$. In view of Proposition 15 and Lemma 4 the Schauder
basis exists in $Y_{K,n_0}$ and consequently, in $M_{\Lambda ,p}$ as
well.


\begin{thebibliography}{199}

\bibitem{ialalam2008} I. Al Alam. "A M\"untz space having no
complement in $L_1$". Proc. Amer. Math. Soc. {\bf 136: 1} (2008),
193-201.

\bibitem{barib} N.K. Bari. "Trigonometric series" (Oxford: Pergamon Press, 1964).

\bibitem{borwerdb} P. Borwein, T. Erd\'elyi. "Polynomials and polynomial
inequalities" (New-York: Springer-Verlag, 1995).

\bibitem{borerdjams1997} P. Borwein, T. Erd\'elyi. "Generalizations
of M\"untz's theorem via a Remez-type inequality for M\"untz
spaces". J. Amer. Mathem. Soc. {\bf 10: 2} (1997), 327-349.

\bibitem{bruijnb} N.G. de Bruijn. "Asymptotic methods in Analysis"
(Amsterdam: North Holland Publishiung Co, 1958).

\bibitem{clarkerd43} J.A. Clarkson, P. Erd\"os. "Approximation by polynomials".
Duke Mathem. J. {\bf 10: 1} (1943), 5-11.

\bibitem{edwardsb} R.E. Edwards. "Functional Analysis. Theory and applications"
(New York: Holt, Rinehart and Winston, 1965).

\bibitem{fichtb} G.M. Fichtenholz. "Differential- und Integralrechnung",
V. 1-3 (Berlin: VEB Deutscher Verlag f\"ur Wissenschaften, 1973).

\bibitem{grinbl41} M.M. Grinblum. "Some theorems on bases in Banach
spaces". Soviet Dokladi {31: 5} (1941), 428-432.

\bibitem{gurlusb} V.I. Gurariy, W. Lusky. "Geometry of M\"untz spaces and related questions".
Lecture Notes in Mathematics, {\bf 1870} (Berlin: Springer, 2005).

\bibitem{gurizv66} V.I. Gurariy. "Bases in spaces of continuous
functions on compacts and some geometrical questions". Math. USSR.
Izvestija {\bf 30: 2} (1966), 289-306.

\bibitem{jarchb} H. Jarchow. "Locally convex spaces" (Stuttgart: B.G. Teubner, 1981).

\bibitem{kolmfomb} A.N. Kolmogorov, S.V. Fomin. "Elements of theory
of functions and functional analysis" (Moscow: Nauka, 1989).

\bibitem{lindliorb} J. Lindenstrauss, L. Tzafriri. "Classical Banach spaces";
V. {\bf 1, 2 }. A series of modern surveys in mathematics {\bf 97}
(Berlin: Springer-Verlag, 1979).

\bibitem{ludksmj} S. V. Ludkovsky. "$\kappa $-normed topological vector spaces".
Siber. Mathem. J. {\bf 41: 1} (2000), 141-154.

\bibitem{ludjms09} S. V. Ludkovsky. "Duality of $\kappa $-normed topological
vector spaces and their applications". J. Mathem. Sci. (New York:
Springer) {\bf 157: 2} (2009), 367-385.

\bibitem{ludlusky15} S.V. Ludkowski, W. Lusky. "On the
geometry of M\"untz spaces". Journal of Function Spaces. online
first, DOI 10.1155/2015/787291.

\bibitem{lusky92} W. Lusky. "On Banach spaces with the commuting bounded approximation property".
Arch. Math. {\bf 58: 6} (1992), 568-574.

\bibitem{lusky96} W. Lusky. "On Banach spaces with bases". J. Funct. Anal.
{\bf 138} (1996), 410-425.

\bibitem{lusky98} W. Lusky. "Three space properties and basis extensions".
Israel. J. Mathem. {\bf 107} (1998), 17-27.

\bibitem{lusky03} W. Lusky. "Three space problems and bounded approximation property".
Stud. Mathem. {\bf 159: 3} (2003), 417-434.

\bibitem{lusky04} W. Lusky. "On Banach spaces with unconditional bases". Israel. J. Math. {\bf 143} (2004),
239-251.

\bibitem{naribeckb} L. Narici, E. Beckenstein. "Topological vector
spaces" (New York: Marcel Dekker, Inc., 1985).

\bibitem{olverb} F.W.J. Olver. "Asymptotics and special functions"
(New York: Academic Press, 1974).

\bibitem{prbrmarb1} A.P. Prudnikov, Yu.A. Brychkov, O.I. Marichev.
"Intergals and series". V.1 (Moscow: Nauka, 1981).

\bibitem{reedsim} M. Reed, B. Simon. "Methods of modern
mathematical physics". V.2 (New York: Academic Press, 1977).

\bibitem{shabatb} B.V. Shabat. "An introduction into complex analysis"
(Moscow: Nauka, 1985).

\bibitem{shirb} A.N. Shiryayev. "Probability" (Moscow: MTzNMO, 2011).

\bibitem{steinb} E.M. Stein. "Singular integrals and
diferentiability properties of functions" (Princeton, NJ: Princeton
University Press, 1986).

\bibitem{schwartzb59} L. Schwartz. "\'Etude des sommes d'exponentielles"; 2-\'eme \'ed. (Paris: Hermann, 1959).

\bibitem{stepanetsb} A.I. Stepanets. "Classification and approximation of periodic
functions", Ser. Mathematics and its applications V. {\bf 333} (Dordrecht: Kluwer Acad. Publ., 1995).

\bibitem{woytb} P. Wojtaszczyk. "Banach spaces for analysts".
Cambridge studies in advanced mathematics, {\bf 25}. (Cambridge: Cambr. Univ. Press,
1991).

\bibitem{zaanb} A.C. Zaanen. "Continuity, integration and Fourier theory" (Berlin: Springer, 1989).

\bibitem{zygmb} A. Zygmund. "Trigonometric series", V. 1, 2, Third Edition
(Cambridge: Cambridge Univ. Press, 2002).

\end{thebibliography}
\end{document}